\documentclass[9pt]{article}

\usepackage{amssymb,latexsym,amsmath}
\usepackage{xcolor}
\usepackage{bm}

\usepackage{hyperref}

\begin{document}

\newtheorem{theorem}{Theorem}[section]

\newtheorem{proposition}[theorem]{Proposition}

\newtheorem{lemma}[theorem]{Lemma}

\newtheorem{corollary}[theorem]{Corollary}

\newtheorem{definition}[theorem]{Definition}

\newtheorem{remark}[theorem]{Remark}

\newtheorem{exempl}{Example}[section]

\newenvironment{exemplu}{\begin{exempl}  \em}{\hfill $\surd$

\end{exempl}}

\newcommand{\ea}{\mbox{{\bf a}}}
\newcommand{\eu}{\mbox{{\bf u}}}
\newcommand{\ep}{\mbox{{\bf p}}}
\newcommand{\ed}{\mbox{{\bf d}}}
\newcommand{\eD}{\mbox{{\bf D}}}
\newcommand{\eK}{\mathbb{K}}
\newcommand{\eL}{\mathbb{L}}
\newcommand{\eB}{\mathbb{B}}
\newcommand{\ueu}{\underline{\eu}}
\newcommand{\ueo}{\overline{u}}
\newcommand{\oeu}{\overline{\eu}}
\newcommand{\ew}{\mbox{{\bf w}}}
\newcommand{\ef}{\mbox{{\bf f}}}
\newcommand{\eF}{\mbox{{\bf F}}}
\newcommand{\eC}{\mbox{{\bf C}}}
\newcommand{\en}{\mbox{{\bf n}}}
\newcommand{\eT}{\mbox{{\bf T}}}
\newcommand{\eV}{\mbox{{\bf V}}}
\newcommand{\eU}{\mbox{{\bf U}}}
\newcommand{\ev}{\mbox{{\bf v}}}
\newcommand{\eve}{\mbox{{\bf e}}}
\newcommand{\uev}{\underline{\ev}}
\newcommand{\eY}{\mbox{{\bf Y}}}
\newcommand{\eP}{\mbox{{\bf P}}}
\newcommand{\eS}{\mbox{{\bf S}}}
\newcommand{\eJ}{\mbox{{\bf J}}}
\newcommand{\leb}{{\cal L}^{n}}
\newcommand{\eI}{{\cal I}}
\newcommand{\eE}{{\cal E}}
\newcommand{\hen}{{\cal H}^{n-1}}
\newcommand{\eBV}{\mbox{{\bf BV}}}
\newcommand{\eA}{\mbox{{\bf A}}}
\newcommand{\eSBV}{\mbox{{\bf SBV}}}
\newcommand{\eBD}{\mbox{{\bf BD}}}
\newcommand{\eSBD}{\mbox{{\bf SBD}}}
\newcommand{\ecs}{\mbox{{\bf X}}}
\newcommand{\eg}{\mbox{{\bf g}}}
\newcommand{\paromega}{\partial \Omega}
\newcommand{\gau}{\Gamma_{u}}
\newcommand{\gaf}{\Gamma_{f}}
\newcommand{\sig}{{\bf \sigma}}
\newcommand{\gac}{\Gamma_{\mbox{{\bf c}}}}
\newcommand{\deu}{\dot{\eu}}
\newcommand{\dueu}{\underline{\deu}}
\newcommand{\dev}{\dot{\ev}}
\newcommand{\duev}{\underline{\dev}}
\newcommand{\weak}{\rightharpoonup}
\newcommand{\weakdown}{\rightharpoondown}
\renewcommand{\contentsname}{ }

\title{Symplectic bipotentials}
\author{\textbf{Mohammad Harakeh}$^1$, \textbf{Michael Ban}$^2$, \textbf{G\'ery de Saxc\'e}$^1$ \\
$^1$ Univ. Lille, CNRS, Centrale Lille, UMR 9013 – LaMcube –\\ Laboratoire de m\'ecanique multiphysique multi\'echelle,\\
France (email: gery.de-saxce@univ-lille.fr)\\
$^2$ Institute of General Mechanics, RWTH Aachen University, Germany}

\maketitle


\begin{abstract}
In a previous paper, we proposed a symplectic version of Brezis-Ekeland-Nayroles principle based on the concepts of Hamiltonian inclusions and symplectic polar functions. We applied it to the standard plasticity. The object of this work is to extend the previous formalism to non associated plasticity. For this aim, we generalize the concept of bipotential to dynamical systems. The keystone idea is to define a symplectic bipotential. We present a method to build it from a bipotential. Next, we generalize the symplectic Brezis-Ekeland-Nayroles principle to non associated dissipative laws. We apply it to the non associated plasticity and to the unilateral contact law with Coulomb's dry friction.
\end{abstract}


\vspace{0.3cm}


\maketitle

\section{Introduction}

Among the tools of the differential geometry, one of the most used is the geometry of the Riemannian spaces, equipped with a symmetric 2-covariant tensor, the metric
$$ g (d\bm{z}, d\bm{z}') = d\bm{z}^T  \bm{K}\, d \bm{z}'
$$
allowing to measure lengths. Another interesting geometry is the one of the symplectic space, equipped with an antisymmetric 2-covariant tensor
$$ \omega (d\bm{z}, d\bm{z}') = d\bm{z}^T  \bm{J}\, d\bm{z}'
$$ 
allowing to measure areas (see \cite{Libermann Marle 1987}, \cite{Souriau 1997b}). One of the origins of the symplectic geometry is the dynamics of reversible systems. We are working with the phase space of which the elements are of the form
$$\bm{z}  =\left[ \begin{array} {c}
                       \bm{x}  \\
                       \bm{y} \\
                    \end {array} \right] 
$$
where $\bm{x}$ are the degrees of freedom and $\bm{y}$ are the corresponding momenta. It is equipped with the canonical symplectic form 
$$\omega (d \bm{z}, d \bm{z}') =
                    \left[ d\bm{x} , d\bm{y} \right]
                    \left[ \begin{array} {cc}
                       \bm{0} & \bm{I}   \\
                     - \bm{I}  & \bm{0}  \\
                    \end {array} \right]\,
                    \left[ \begin{array} {c}
                       d \bm{x}'  \\
                       d \bm{y}' \\
                    \end {array} \right] 
$$
where  $\bm{I}$ is the identity operator.  Introducing the Hamiltonian vector field (or symplectic gradient)
$$ \dot{\bm{z}} = \bm{X}_H = \bm{J} \cdot \nabla_{\bm{z}}  H (t, z )
$$
allows to recover the canonical equations
$$ \dot{\bm{x} }   =   \nabla_{\bm{y} } H, \qquad 
 \dot{\bm{y} } = - \nabla_{\bm{x} } H  .   
$$
Considering the Lagrangian $L$ associated to the Hamiltonian $H$, this formulation can be recast in integral form thanks to

\textbf{Hamilton's variational principle:} \textit{the natural evolution of the systems minimizes}
 \textit{the action} 
  $$\alpha \left[\bm{x}\right]= \int^{T}_0 L (\bm{x}, \dot{\bm{x}}, t) \,\mbox{ d} t $$  
\textit{among the paths} $t \mapsto \bm{x} (t)$ \textit{satisfying}
$$\bm{x} (0) = \bm{x}_0, \quad 
    \bm{x} (T) = \bm{x}_T 
$$ 
It is important to remark that the application of this principle is strictly limited to the reversible systems but unfortunately \textit{it fails} for dissipative ones. 

For dissipative systems, several authors proposed unified frameworks inspired from both  previously mentioned geometries. We quickly review now these theoretical frameworks. To avoid indulging in very complex details, we present them in the smooth case (although we shall consider later on much more general extensions to the non-smooth mechanics):

\begin{itemize}
\item The {\it metriplectic systems} were introduced by Morrison (see \cite{Barbaresco 2022}, \cite{Coquinot 2020}, \cite{Materassi 2017}, \cite{Morrison 1986}) and developed further by Grmela and \"Ottinger (see \cite{Grmela 1997}, \cite{Ottinger 1997}) as the {\it GENERIC systems} (General Equation for Non-Equilibrium Reversible-Irreversible Coupling). They combine the Hamiltonian formulation and the Onsager one, according to the evolution law
$$ \dot{\bm{z}} = \bm{J} \cdot \nabla_{\bm{z} } H (\bm{z} ) 
                        + \bm{K}   \cdot \nabla_{\bm{z} } S (\bm{z})
$$
where Onsager term is built from a symmetric and  positive-definite operator $\bm{K}$ and an
entropy-like function $S$. A variational formulation of GENERIC can be found in \cite{Manh Hong Duong 2013}.
\item The {\it Port-Hamiltonian systems} were introduced by Brockett \cite{Brockett 1977} and van der Schaft \cite{van der Schaft 1984}
$$ \dot{\bm{z}} = (\bm{J} - \bm{R})\cdot \nabla_{\bm{z}} H (\bm{z} )
$$  
where the symmetric and positive-definite operator $ \bm{R}$ modelizes the resistive effects (because of the minus sign).                   
\item The {\it rate-independent systems} proposed by Mielke and Theil \cite{mielketh99}, Mielke \cite{mielke} and developped with applications in Mielke and Roub\'{\i}\v{c}ek \cite{MR06b}, are based on two fundamental conditions:
\begin{itemize}
\item The stability condition: 
$$ \nabla_{\bm{x}} E (\bm{x}) \cdot \bm{w} + \Phi (\bm{w}) \geq  0, \; \forall \bm{w} $$ 
\item The power balance: 
$$ \nabla_{\bm{x}} E (\bm{x}) \cdot \dot{\bm{x}} + \Phi (\dot{\bm{x}} ) = 0 $$ 
\end{itemize}
where $ E $ is the energy functional and $\Phi$ is a 1-homogeneous dissipation potential depending on the velocity.
\item The {\it Hamiltonian inclusions}  were proposed by Buliga \cite{bham},
\begin{equation}
 \dot{\bm{z}} = \bm{J} \cdot \nabla_{\bm{z} } H (\bm{z} ) +  \bm{J} \cdot \nabla_{\dot{\bm{z}} } \Phi (\dot{\bm{z}} ) 
\label{dot(z) = J nabla_z H (z) + nabla_(dot(z)) Phi (dot(z))}
\end{equation}
where $ \Phi$ is a convex dissipation potential. Because of the first term, It is clearly related to the two former formalisms (metriplectic and Port-Hamiltonian systems), but it is also inspired from the rate-independent systems because of the dissipation potential. Nevertheless, it is important to remark that $\Phi$ is not necessarily 1-homogeneous or even homogeneous. 
\end{itemize}

On the other hand, the Brezis-Ekeland-Nayroles variational principle (BEN in short) was proposed independently one of each other by Brezis and Ekeland \cite{Brezis Ekeland 1976}, and Nayroles \cite{Nayroles 1976}. 
In \cite{Stefanelli 2008}, Stefanelli  used the BEN principle to represent the quasistatic evolution of an elastoplastic material with hardening in order to prove the convergence of time and space-time discretizations as well as to provide some possible a posteriori error control. Finally, Ghoussoub and MacCann characterized the path of steepest descent of a non-convex potential as the global minimum of BEN functional \cite{Ghoussoub 2004}.

Latter on, Buliga and de Saxc\'e merged, in a symplectic framework, the formalism of Hamiltonian inclusions with BEN principle to extend it to dynamics \cite{SBEN}. This generalized approach is quite naturally called symplectic BEN principle (SBEN in short).
The key idea is to decompose additively the time rate $\dot{\bm{z}}$ into reversible part $\dot{\bm{z}}_R$ (the symplectic gradient) and dissipative or irreversible one $\dot{\bm{z}}_I$, and then to define the symplectic subdifferential $\partial^{\omega} \Phi (\bm{z})$ of the dissipation potential. To release the restrictive hypothesis of $1$-homogeneity (in particular to address viscoplasticity), we introduce in this work the symplectic Fenchel polar $\phi^{*\omega}$, that allows to  build theoretical methods to model and analyze dynamical dissipative systems in a consistent geometrical framework with the numerical approaches not very far in the background. Numerical simulations  with the BEN principle were performed for elastoplastic structures in statics (\cite{Cao 2020}, \cite{Cao 2021a}) and in dynamics \cite{Cao 2021b}.

An advantage of BEN principle is the easiness to be generalized. Indeed, it is worth to know that many realistic dissipative laws, called non-associated, cannot be cast in the mould of the standard ones deriving \textcolor{black}{from} a dissipation potential. To skirt this pitfall, the author proposed in \cite{saxfeng} a new theory based on a function called bipotential. It represents physically the dissipation and generalizes the sum of the dissipation potential and its Fenchel polar, reason for which extension of BEN principle is natural. 
The applications of the bipotential approach to solid Mechanics are various: Coulomb's friction law \cite{sax CRAS 92}, non-associated Drucker-Prager  \cite{sax boussh IJMS 98} and Cam-Clay models \cite{Zouain 2010} in Soil Mechanics, cyclic Plasticity (\cite{bodo sax EJM 01}, \cite{Bouby 2009}, \cite{Bouby 2015}, \cite{sax CRAS 92}, \cite{Magnier 2006}) and Viscoplasticity \cite{hjiaj bodo CRAS 00} of metals with non linear kinematical hardening rule, Lemaitre's damage law \cite{bodo},  the coaxial laws (\cite{dangsax}, \cite{vall leri CONST 05}). Such kind of materials are called implicit standard materials. A synthetic review of these laws can be found in the two latter references. It is also worth to notice that monotone laws that \textcolor{black}{do} not admit a convex potential can be represented by Fitzpatrick's function \cite{Fitzpatrick 1988} which is a bipotential.

The paper is organized as follows. Section 'Symplectic subdifferential and symplectic polar function' is a reminder of the concepts of symplectic subdifferential and symplectic polar function, with attention paid to the particular case of interest when the duality is a scalar product. Section 'Hamiltonian inclusions and symplectic BEN principle' is a short presentation of the formalism of Hamiltonian inclusions and the symplectic BEN principle. In Section 'Symplectic bipotentials', we recall the notion of bipotential and we extend it to the dynamics by defining the symplectic bipotential. We propose a method of  construction of the symplectic bipotential from a bipotential in the particular case of interest mentioned above. In Section 'BEN principle for symplectic bipotentials', we state the BEN principle for symplectic bipotentials and we apply it to the non associated plasticity and to the crack extension with friction contact between the crack sides.

\section{Symplectic subdifferential and symplectic polar function}
\label{Section - Symplectic subdifferential and symplectic polar function}

\subsection{Basic definitions}

To begin with, we recall the basic definitions of \cite{bham} and \cite{SBEN}. $X$ and $Y$ are topological, locally convex, real vector spaces of dual variables $\bm{x} \in X$ and $\bm{y} \in Y$ through the dual pairing
$$\langle \cdot , \cdot \rangle : X \times Y \rightarrow \mathbb{R}
$$
such that  any  continuous linear functional on $X$ (resp. on $Y$) has the form $\bm{x}  \mapsto \langle \bm{x} ,\bm{y} \rangle$, for some $\bm{y}  \in Y$ (resp.  $\bm{y}  \mapsto \langle \bm{x} ,\bm{y} \rangle$, for some  $\bm{x}  \in X$)\footnote{For instance if $X$ and $Y$ are Hilbert spaces}. Usually, $\bm{x}$ are the degrees of freedom (that may be fields) and $\bm{y}$ are the corresponding dynamic momenta. A generic element of $X \times Y$ is denoted $\bm{z}  = (\bm{x} ,\bm{y} )$. The space $X \times Y$ is equipped with a natural symplectic form 
$\omega: \left( X \times Y \right)^{2} \rightarrow \mathbb{R}$ which is defined via the duality product for any $\bm{z} =(\bm{x} ,\bm{y} )$ and  $\bm{z} '=(\bm{x}',\bm{y}')$ by   
$$\omega ( \bm{z} , \bm{z}' ) \, = \langle \bm{x} , \bm{y}' \rangle - \langle \bm{x} ' , \bm{y}  \rangle 
$$
Let $F: X \times Y \rightarrow \mathbb{R}\cup \left\{+\infty\right\}$ be a convex lower semicontinuous  function.  The symplectic subdifferential of $F$ at $\bm{z}$ where $F$ has a finite value is the set
\begin{eqnarray}
    \partial^{\omega} F (\bm{z})  =  
\{ 
\bm{z}' \in X \times Y \mbox{ : } \forall  \bm{z}'' \in X \times Y \quad \nonumber \\ 
F(\bm{z}+\bm{z}'') - F(\bm{z}) \geq
\omega (\bm{z}' , \bm{z}'') \} 
\label{dssub}
\end{eqnarray}
The elements of the symplectic subdifferential are called symplectic subgradients and satisfy the so-called Hamiltonian inclusion
\begin{equation}
   \bm{z}' \in \partial^{\omega} F (\bm{z}) 
   \label{z' in partial^omega F (z)}
\end{equation}
The symplectic polar of $F$ is the function 
\begin{equation}
F^{*\omega}(\bm{z}') =  \sup \left\{ \omega(\bm{z}',\bm{z}) - F(\bm{z}) \mbox{ : } \bm{z} \in X \times Y \right\}
\label{dspolar}
\end{equation}
By construction of the symplectic polar, the symplectic Fenchel inequality holds
\begin{equation}
\forall  \bm{z}, \bm{z}' \in X \times Y, \qquad
  F(\bm{z}) + F^{*\omega}(\bm{z}') - \omega(\bm{z}',\bm{z}) \geq 0
\label{symplectic Fenchel inequality}
\end{equation}
The equality in the previous relation is reached when $\bm{z}'$ is a symplectic subgradient of $F$ at $\bm{z}$. Indeed, it has been proved  in \cite{SBEN} that 
\begin{equation}
  F(\bm{z}) + F^{*\omega}(\bm{z}') - \omega(\bm{z}',\bm{z}) = 0
  \quad \Leftrightarrow \quad 
 \bm{z}' \in \partial^{\omega} F (\bm{z})
\label{extremality <-> dot(z)_I in partial^omega F(dot(z))}
\end{equation} 
In such case, we say that the couple is extremal.  
We would like now to find the inverse law of (\ref{z' in partial^omega F (z)}) in term of the symplectic polar. Because of the symplectic Fenchel inequality (\ref{symplectic Fenchel inequality}), we have
$$  F(\bm{z}) + F^{*\omega} (\bm{z}' + \bm{z}'')  \geq \omega (\bm{z}' +\bm{z}'', \bm{z}) 
$$
If $\bm{z}'$ is a symplectic subgradient of $F$ at $\bm{z}$, (\ref{extremality <-> dot(z)_I in partial^omega F(dot(z))}) gives
$$  F(\bm{z}) + F^{*\omega}(\bm{z}') - \omega(\bm{z}',\bm{z}) = 0
$$
By substraction of the two previous relations member by member, we have
$$  F^{*\omega} (\bm{z}' + \bm{z}'') - F^{*\omega} (\bm{z}')  \geq \omega (- \bm{z}, \bm{z}'') 
$$
Taking into account the definition of the symplectic subdifferential, we obtain the inverse law of (\ref{z' in partial^omega F (z)})
\begin{equation}
   - \bm{z} \in \partial^{\omega} F^{*\omega} (\bm{z}') 
   \label{- z in partial^omega F^(*omega) (z')}
\end{equation}
By similar arguments to those developed in (\cite{SBEN}), we can prove the equivalence of this differential inclusion with (\ref{extremality <-> dot(z)_I in partial^omega F(dot(z))}), then it holds
\begin{eqnarray}
  F(\bm{z}) + F^{*\omega}(\bm{z}') - \omega(\bm{z}',\bm{z}) = 0
  \quad \Leftrightarrow \quad \nonumber \\
 \bm{z}' \in \partial^{\omega} F (\bm{z})
  \quad \Leftrightarrow \quad 
 - \bm{z} \in \partial^{\omega} F^{*\omega} (\bm{z}')
\label{extremality <-> dot(z)_I in partial^omega F(dot(z)) <-> dot(z) in partial^omega F^(*omega) (dot(z))}
\end{eqnarray} 

\subsection{Particular case of interest}
\label{SubSection Particular case of interest}

In the applications, we shall be led to consider the simpler case where  $X=Y$ is a Hilbert space, then the duality  $\displaystyle \langle \cdot, \cdot \rangle$ is a scalar product.  Moreover the space $X \times Y$ is dual with itself, with the scalar product: 
\begin{equation}
    \langle \langle (\bm{x},\bm{y}), (\bm{x}', \bm{y}') \rangle \rangle 
    \, = \, \langle \bm{x}, \bm{x}' \rangle + \langle \bm{y}, \bm{y}' \rangle 
\label{scalar product << . , . >>}
\end{equation}
We  introduce now a continuous endomorphism of $X \times Y$
$$\bm{J}: X \times Y \rightarrow X \times Y: (\bm{x},\bm{y}) \mapsto \bm{J} (\bm{x},\bm{y}) = (\bm{y}, - \bm{x}).
$$ 
The map $\bm{J}$ has the properties:
\begin{equation}
  \bm{J}^2 = - \bm{I}_{X \times Y} , \qquad 
  \bm{J}^{-1} =  \bm{J}^T = - \bm{J} 
\label{J^2 = - I & J^(-1) = J^T = - J}
\end{equation}
and
\begin{equation}
   \omega(\bm{z}, \bm{z}') \, = \, \langle \langle \bm{z}, \bm{J} \, \bm{z}' \rangle \rangle 
                                             = \, \langle \bm{x}, \bm{y}' \rangle - \langle \bm{x}', \bm{y} \rangle 
\label{omega(z, z') = << z, J z' >>}
\end{equation}

The subdifferential of $F$ is by definition: 
\begin{eqnarray}
    \partial F (\bm{z})  =  
\{ 
\bm{z}' \in X \times Y \mbox{ : } \forall \, \bm{z}'' \in X \times Y \quad \nonumber \\ 
F(\bm{z}+\bm{z}'') - F(\bm{z}) \geq
\langle \langle\bm{z}' , \bm{z}'' \rangle \rangle \}  \nonumber
\end{eqnarray}
By comparison with the definition (\ref{dssub}) of the symplectic subdifferential of $F$, we obtain in the particular case $X = Y$
\begin{equation}
    \bm{z}' \in \partial^{\omega} F(\bm{z}) \quad \Leftrightarrow \quad     \bm{z}' \in \bm{J} \, \partial F(\bm{z}) 
\label{z' in partial^omega F (z) <-> z' in J partial F (z)}
\end{equation}
$$\bm{z}' \in   \partial F(\bm{z}) \quad \Leftrightarrow \quad - \bm{z}' \in \bm{J} \,  \partial^{\omega} F(\bm{z})
$$
Likewise, we recall the definition  of the Fenchel polar of $F$
$$F^{*}(\bm{z})  =  \sup \left\{ \langle \langle \bm{z}', \bm{z} \rangle \rangle  -  F(\bm{z}') \mbox{ : } \bm{z}' \in X \times Y \right\} 
$$ 
Remark that, from definition (\ref{dspolar}) of the symplectic polar, we get in the particular case $X =Y$
\begin{equation}
    F^{*\omega}  =  F^{*} \circ  \bm{J}^{-1}, \qquad
   F^{*} = F^{*\omega}   \circ \bm{J}
\label{F^(*omega) = F^* circ J^(-1) & F^* = F^(*omega) circ J}
\end{equation}

We can verify that 
$$ - \bm{z} \in \partial^{\omega} F^{*\omega} (\bm{z}') 
     \quad \Leftrightarrow \quad 
     \bm{z} \in \partial F^* (\bm{J}^{-1} \, \bm{z}') 
$$
Indeed, considering a continuous homomorphism $\bm{A}$ of vector spaces, we use the chain rule for subdifferentials (see for instance \cite{Ekeland Temam 1999}, page 27, Proposition 5.7)
\begin{equation}
     \partial (F \circ \bm{A}) (\bm{z}) = \bm{A}^T \partial \, F (\bm{A} \, \bm{z})
\label{chain rule}
\end{equation}
Owing to (\ref{z' in partial^omega F (z) <-> z' in J partial F (z)}), (\ref{F^(*omega) = F^* circ J^(-1) & F^* = F^(*omega) circ J}) and
 applying the chain rule to $ \bm{A} = \bm{J}^{-1}$, one has
$$ \partial^{\omega} F^{*\omega} (\bm{z}')  = \bm{J} \, \partial (F^{*} \circ  \bm{J}^{-1} ) (\bm{z}') 
   =  \bm{J} \, (\bm{J}^{-1})^T \partial F^{*} ( \bm{J}^{-1} \bm{z}') 
$$
in which, taking into account (\ref{J^2 = - I & J^(-1) = J^T = - J})
\begin{equation}
   \bm{J} \, (\bm{J}^{-1})^T = \bm{J} \, (\bm{J}^T )^T = \bm{J}^2 = - \bm{I}_{X \times Y}
\label{J (J^-1)^T = - I_(X times Y)}
\end{equation}
then
$$ \partial^{\omega} F^{*\omega} (\bm{z}')  =  - \partial F^{*} ( \bm{J}^{-1} \bm{z}') 
$$

\section{Hamiltonian inclusions and symplectic BEN principle}

\label{Section - Hamiltonian inclusions and symplectic BEN principle}

Coming back to the general setting, we do not suppose that $X = Y$, then the operator $\bm{J}$ no longer makes sense. 
For the evolution of a dissipative dynamical system, we suppose an additive decomposition of the velocity in reversible and irreversible parts: 
$$ \dot{\bm{z}} =  \dot{\bm{z}}_{R} + \dot{\bm{z}}_{I},\quad 
\dot{\bm{z}}_{R} = \bm{X}_H,\quad 
\dot{\bm{z}}_{I}  =  \dot{\bm{z}}  - \bm{X}_H 
$$
where the reversible part of the velocity is given by the symplectic gradient of the Hamiltonian $H$ (or Hamiltonian vector field)
$$ \omega (\bm{X}_H, \bm{\zeta}) = \lim_{\epsilon \to 0} \frac{1}{\epsilon} 
      (H (t, \bm{z} + \epsilon \, \bm{\zeta}) 
      - H (t, \bm{z}) )
$$
If $\bm{z}$ is a field, the symplectic gradient is a functional derivative. Besides, we suppose that the dissipative flow rule is an Hamiltonian  inclusion
\begin{equation}
\dot{\bm{z}}_{I} \in \partial^\omega \phi (\dot{\bm{z}})
\label{dot(z)_I in partial^omega F(dot(z))}
\end{equation}
deriving from a disipation potential $\phi$. Hence the path $\left[0, T\right] \rightarrow X \times Y: t \mapsto \bm{z} (t)$ representing the behavior of the system is governed by the following constitutive law
$$ \dot{\bm{z}} \in \bm{X}_H + \partial^\omega \phi (\dot{\bm{z}})
$$
with initial conditions $\bm{z} (0) = \bm{z}_0$ and possible supplemental conditions (for instance, boundary conditions if it is a field), in which case it is called an admissible path. Remarking that the natural path must satisfy the Hamiltonian inclusion (\ref{dot(z)_I in partial^omega F(dot(z))}) almost everywhere in $\left[0, T\right]$ and taking into account the symplectic Fenchel inequality (\ref{symplectic Fenchel inequality}) and the equivalence  (\ref{extremality <-> dot(z)_I in partial^omega F(dot(z))}), the symplectic Brezis-Ekeland-Nayroles principle (in short SBEN principle) was stated in \cite{SBEN}:

\textbf{SBEN principle.} \textit{The natural evolution of the system minimizes the functional}
\begin{equation}
\Pi(\bm{z}) = \int_{0}^{T} \left\{ \phi(\dot{\bm{z}})      
    + \phi^{*\omega}(\dot{\bm{z}} - \bm{X}_H) 
    - \omega (\dot{\bm{z}} - \bm{X}_H, \dot{\bm{z}})  \right\} \mbox{ dt} 
\label{SBEN 1 Pi (z) =}
\end{equation}
\textit{among the admissible paths $t \mapsto \bm{z} (t) $ and the minimum is zero.} 

This principle is an extension of the original BEN principle  (\cite{Brezis Ekeland 1976} , \cite{Nayroles 1976})  to the dynamics.

\section{Symplectic bipotentials}

\label{Section - Symplectic bipotentials}

The bipotential is a tool proposed in \cite{sax CRAS 92} to modelize the non associated dissipative laws. A rigorous mathematical formalism was proposed in (\cite{Buliga 2008}, \cite{Buliga 2010}). A bipotential is a function $b: \mathcal{X} \times \mathcal{Y} \rightarrow \mathbb{R} \cup \left\lbrace +\infty \right\rbrace$ such that:
\begin{itemize}
    \item [(i)] $b$ is bi-convex (\textit{i.e.} convex with respect to each of its arguments) and bi-lower semicontinuous
    \item  [(ii)]   $\forall \bm{x}' \in \mathcal{X}, \bm{y}' \in \mathcal{Y},\qquad b(\bm{x}', \bm{y}') \geq \left\langle \bm{x}', \bm{y}' \right\rangle$
    \item   [(iii)]   If one of the following relations is true, the other ones also are
    \begin{eqnarray}
        b(\bm{x}, \bm{y}) = \left\langle \bm{x}, \bm{y} \right\rangle 
                \qquad \Leftrightarrow \qquad  \nonumber \\ 
                \bm{y} \in \partial b (\bullet, \bm{y}) (\bm{x})
                \qquad \Leftrightarrow \qquad  
                \bm{x} \in \partial b (\bm{x}, \bullet) (\bm{y})  \nonumber
    \end{eqnarray}
   The meaning of the second relation is that if we differentiate with respect to $\bm{x}$ at constant $\bm{y}$, then $\bm{y}$ is an implicit   function of $\bm{x}$. 
\end{itemize}
In the particular event of associated laws, the bipotential is separated
$$ b(\bm{x}, \bm{y}) = \phi(\bm{x}) + \phi^* ( \bm{y})
$$
where the dissipation potential $\phi$ and its Fenchel polar $\phi^*$ are convex and lower semicontinuous, the equivalence of the three relations of (iii) is true in general as a consequence of (i) and (ii). It is not so for a non separated bipotential, reason for which it is an axiom of the definition of bipotentials. 

For the non associated laws in the context of dynamics, we propose to define a symplectic bipotential as a function $\hat{b}: (X \times Y)^2 \mapsto \mathbb{R} \cup \left\lbrace +\infty \right\rbrace$ such that:
\begin{itemize}
    \item [(a)] $\hat{b} $ is bi-convex and bi-lower semicontinuous
    \item [(b)] $\forall \dot{\bm{z}}', \dot{\bm{z}}'_I \in Z,\qquad \hat{b} (\dot{\bm{z}}'_I, \dot{\bm{z}}' )                                                      \geq \omega (\dot{\bm{z}}'_I, \dot{\bm{z}}' )$
    \item [(c)] \begin{eqnarray}
     \hat{b} (\dot{\bm{z}}_I, \dot{\bm{z}}) = \omega (\dot{\bm{z}}_I, \dot{\bm{z}} )
                \quad \Leftrightarrow \qquad\qquad \nonumber \\ 
                \dot{\bm{z}}_I \in \partial^{\omega} \hat{b}  (\dot{\bm{z}}_I, \bullet) (\dot{\bm{z}})
                \quad \Leftrightarrow \quad  
                - \dot{\bm{z}} \in \partial^{\omega} \hat{b}  (\bullet, \dot{\bm{z}}) (\dot{\bm{z}}_I)   \nonumber  
    \end{eqnarray}
    \end{itemize}
In the particular case of interest $X=Y$ with $X \times Y$ endowed with the scalar product $\langle \langle \cdot , \cdot \rangle \rangle$, we can verify that if $b: (X \times Y)^2 \rightarrow \mathbb{R}\cup \left\{+\infty\right\} $ is a bipotential 
\begin{equation}
   \hat{b} (\dot{\bm{z}}_I, \dot{\bm{z}}) = b (\bm{J}^{-1} \dot{\bm{z}}_I, \dot{\bm{z}})
\label{hat(b) (dot(z)_I, dot(z)) = b (J^(-1) dot(z)_I, dot(z)) }
\end{equation}
is a symplectic bipotential. 

Indeed, (i) applied to $b$ entails (a) for $\hat{b}$ because the map $\bm{J}$ is linear.

Moreover, (ii) entails (b) if we remark that
\begin{equation}
   \omega (\dot{\bm{z}}_I, \dot{\bm{z}} ) 
   = \langle \langle \dot{\bm{z}}_I , \bm{J}\, \dot{\bm{z}} \rangle \rangle   
   = \langle \langle \bm{J}^T \, \dot{\bm{z}}_I ,  \dot{\bm{z}} \rangle \rangle
   = \langle \langle \bm{J}^{-1} \, \dot{\bm{z}}_I ,  \dot{\bm{z}} \rangle \rangle
\label{omega (dot(z)_I, dot(z)) = << J^(-1) dot(z)_I, dot(z) >>}
\end{equation}

Finally, we have to prove the equivalence between the three following relations
\begin{equation}
    \hat{b} (\dot{\bm{z}}_I, \dot{\bm{z}}) = \omega (\dot{\bm{z}}_I, \dot{\bm{z}} )
\label{hat(b) (dot(z)_I, dot(z)) = omega (dot(z)_I, dot(z))}
\end{equation}
\begin{equation}
               \dot{\bm{z}}_I \in \partial^{\omega} \hat{b}  (\dot{\bm{z}}_I, \bullet) (\dot{\bm{z}})
\label{dot(z)_I in partial^omega hat(b) (dot(z)_I, . ) (dot(z))}
\end{equation}
\begin{equation}
         - \dot{\bm{z}} \in \partial^{\omega} \hat{b}  (\bullet, \dot{\bm{z}}) (\dot{\bm{z}}_I) 
\label{- dot(z) in partial^omega hat(b) ( . , dot(z)) (dot(z)_I)}
\end{equation}

Firstly, owing to (\ref{hat(b) (dot(z)_I, dot(z)) = b (J^(-1) dot(z)_I, dot(z)) }) and (\ref{omega (dot(z)_I, dot(z)) = << J^(-1) dot(z)_I, dot(z) >>}), (\ref{hat(b) (dot(z)_I, dot(z)) = omega (dot(z)_I, dot(z))}) is equivalent to
           \begin{equation}
    b (\bm{J}^{-1} \dot{\bm{z}}_I, \dot{\bm{z}}) = \langle \langle \bm{J}^{-1} \, \dot{\bm{z}}_I ,  \dot{\bm{z}} \rangle \rangle
\label{b (J^(-1) dot(z)_I, dot(z)) = << J^(-1) dot(z)_I, dot(z) >>}
\end{equation}    
Applying (iii), this relation is equivalent to
$$ \bm{J}^{-1} \dot{\bm{z}}_I \in \partial \,  b \,  (\bm{J}^{-1} \dot{\bm{z}}_I, \bullet) (\dot{\bm{z}} )
$$
or,  taking into account (\ref{z' in partial^omega F (z) <-> z' in J partial F (z)}) and (\ref{hat(b) (dot(z)_I, dot(z)) = b (J^(-1) dot(z)_I, dot(z)) })
$$  \dot{\bm{z}}_I \in \bm{J} \, \partial \,  b \,  (\bm{J}^{-1} \dot{\bm{z}}_I, \bullet) (\dot{\bm{z}} )
        = \partial^{\omega}  b \,  (\bm{J}^{-1} \dot{\bm{z}}_I, \bullet) (\dot{\bm{z}} )
$$
$$  \dot{\bm{z}}_I \in 
         \partial^{\omega} \hat{b}  (\dot{\bm{z}}_I, \bullet) (\dot{\bm{z}})
$$
that proves the equivalence between (\ref{hat(b) (dot(z)_I, dot(z)) = omega (dot(z)_I, dot(z))}) and (\ref{dot(z)_I in partial^omega hat(b) (dot(z)_I, . ) (dot(z))}).

Secondly, (\ref{hat(b) (dot(z)_I, dot(z)) = omega (dot(z)_I, dot(z))}) is equivalent to (\ref{b (J^(-1) dot(z)_I, dot(z)) = << J^(-1) dot(z)_I, dot(z) >>}) then, applying (iii), equivalent to 
$$ \dot{\bm{z}}  \in \partial \, b \, ( \bullet, \dot{\bm{z}} ) (\bm{J}^{-1} \dot{\bm{z}}_I)
$$
Besides, owing to (\ref{J (J^-1)^T = - I_(X times Y)}) and the chain rule (\ref{chain rule}), one has
$$  \partial \, b \, ( \bullet, \dot{\bm{z}} ) (\bm{J}^{-1} \dot{\bm{z}}_I)
     = -  \bm{J} \, (\bm{J}^{-1})^T \partial \,  b \,  (\bullet, \dot{\bm{z}} ) ( \bm{J}^{-1} \bm{z}_I) 
$$
$$  \partial \, b \, ( \bullet, \dot{\bm{z}} ) (\bm{J}^{-1} \dot{\bm{z}}_I)
     = -  \bm{J} \, \partial \,  \hat{b} \,  (\bullet, \dot{\bm{z}} ) ( \bm{z}_I)
$$
and, using (\ref{z' in partial^omega F (z) <-> z' in J partial F (z)})
$$    \partial \, b \, ( \bullet, \dot{\bm{z}} ) (\bm{J}^{-1} \dot{\bm{z}}_I)    
     = - \partial^{\omega} \hat{b}  (\bullet, \dot{\bm{z}}) (\dot{\bm{z}}_I) 
$$
that proves the equivalence between (\ref{hat(b) (dot(z)_I, dot(z)) = omega (dot(z)_I, dot(z))}) and (\ref{- dot(z) in partial^omega hat(b) ( . , dot(z)) (dot(z)_I)}).

The relation (\ref{hat(b) (dot(z)_I, dot(z)) = b (J^(-1) dot(z)_I, dot(z)) }) provides a convenient way to construct symplectic bipotentials in dynamical conditions from bipotentials used in statics. 

\section{BEN principle for symplectic bipotentials}

\label{Section - BEN principle for symplectic bipotentials}

Now we come back to the calculus of variations.  The domain of application of the SBEN principle in Mechanics is very broad but limited to problems where the dissipative laws of materials are generated by a potential. For atypical  dissipative laws called non associated, there is no corresponding potential  but a great number of them can be generated by a bipotential. In (\cite{Visintin 2008}, \cite{Visintin 2013}), Visintin extended the Brezis-Ekeland-Nayroles to Fitzpatrick's function, a bipotential for maximal monotone laws.  Our aim is to generalize the SBEN principle to such constitutive laws. The  key idea is to replace in the functional the sum of the potential of dissipation and its symplectic polar by a symplectic bipotential $\hat{B} $. 
 
\textbf{SBEN principle for symplectic bipotentials.} \textit{The natural evolution of the system minimizes the functional}
\begin{equation}
\Pi(\bm{z}) = \int_{0}^{T} \left\{ \hat{B}       
    (\dot{\bm{z}} - \bm{X}_H , \dot{\bm{z}})
    - \omega (\dot{\bm{z}} - \bm{X}_H, \dot{\bm{z}})  \right\} \mbox{ dt} 
\label{SBEN 2 Pi (z) =}
\end{equation}
\textit{among the admissible paths $t \mapsto \bm{z} (t) $ and the minimum is zero.} 

To illustrate the application of this variational principle, we chose two simple situations involving non associated laws.

\vspace{0.5cm}

\subsection{Application to the non associated plasticity}
\label{SubSection Application to the non associated plasticity}

Non associative flow rules in plasticity of geomaterials are ubiquitous. In previous works, we showed that many of them can be generalized by a bipotential. Our goal is to study the plastic flow evolution in the structure during a time interval $[0,T]$.
As usual under the  hypothesis of small strains, we decompose the strain tensor additively in elastic and plastic parts
$$ \bm{\varepsilon} = \bm{\varepsilon}^e + \bm{\varepsilon}^p
$$
Let $\Omega $ be a bounded, open set of $\mathbb{R}^3$ with piecewise smooth boundary $\partial \Omega $, divided into two disjoint parts, $\Gamma_0 $ (called support) where the displacements are imposed, and $\Gamma_1 $ where the surface forces are imposed. The elements of the space $X$ are the couples $\bm{x} =  (\bm{u}, \bm{\varepsilon}^p) $ of fields of displacement and plastic strain. The corresponding dynamical momenta are denoted $\bm{y} = (\bm{p} , \bm{\pi}) $. The dual pairing is of the form
$$ \left\langle \dot{\bm{x}} , \dot{\bm{y}} \right\rangle 
 = \left\langle \dot{\bm{u}} , \dot{\bm{p}} \right\rangle + \left\langle \dot{\bm{\varepsilon}}^p , \dot{\bm{\pi}} \right\rangle 
= \int_\Omega \dot{\bm{u}} \cdot \dot{\bm{p}}  \, d\Omega
  + \int_\Omega \dot{\bm{\varepsilon}}^p : \dot{\bm{\pi}}  \, d\Omega
$$
The Hamiltonian of the structure is 
\begin{eqnarray}
    H (t, \bm{z}) 
    = \int_{\Omega} \left\lbrace \dfrac{1}{2 \rho} \parallel \bm{p} \parallel^ 2 
    + w\, (\nabla \bm{u} - \bm{\varepsilon}^p) - \bm{f} (t)\cdot\bm{u} \right\rbrace  \, d\Omega \nonumber \\
              - \int_{\Gamma_1} \bar{\bm{f}} (t)\cdot \bm{u}  \, d\Gamma \nonumber
\end{eqnarray}
The first term is the kinetic energy, $w$ is  the quadratic elastic strain energy expressed in term of the compliance matrix $\bm{S}$
$$ w (\bm{\varepsilon}) =   \frac{1}{2} \, \bm{\varepsilon} : \bm{S}^{-1} \bm{\varepsilon}
$$
$\bm{f}$ is the volume force and $\bar{\bm{f}}$ is the surface force on the part $\Gamma_1$ of the boundary, the displacement field being equal to an imposed value $\bar{\bm{u}}$ on the remaining part $\Gamma_0$. The stress field is given by the elasticity law
\begin{equation}
\label{elasticity law}
      \bm{\sigma} = \frac{\partial w}{\partial \bm{\varepsilon}} \,  (\nabla \bm{u} - \bm{\varepsilon}^p)
                          = \bm{S}^{-1} \, : (\nabla \bm{u} -\bm{\varepsilon}^p )
\end{equation}
We consider the canonical symplectic form
$$ \omega \, (\dot{\bm{z}}, \dot{\bm{z}}') =  \omega \, ( (\dot{\bm{x}}, \dot{\bm{y}}) , (\dot{\bm{x}}', \dot{\bm{y}}'))
    = \left\langle  \dot{\bm{x}}, \dot{\bm{y}}'  \right\rangle - \left\langle \dot{ \bm{x}}', \dot{\bm{y}} \right\rangle
$$
Denoting the functional derivative by $D$, the symplectic gradient of the Hamiltonian is
$$ \bm{X}_H = ((D_{\bm{p}} H,D_{\bm{\pi}} H), (- D_{\bm{u}} H, - D_{\bm{\varepsilon}^p} H))
$$
$$ \bm{X}_H = ( ( \bm{p} / \rho, \bm{0} ) , ( \nabla \cdot \bm{u} + \bm{f} , \bm{\sigma} ) )
$$

Thus one has
$$ \dot{\bm{z}}_I = \dot{\bm{z}} - \bm{X}_H = ( ( \bm{u}_I , (\bm{\varepsilon}^p)_I ) , ( \bm{p}_I, \bm{\pi}_I ) )
$$
with
$$ \dot{\bm{u}}_I =\dot{\bm{u}} -  \bm{p} / \rho , \quad
     (\dot{\bm{\varepsilon}}^p)_I = \dot{\bm{\varepsilon}}^p , 
$$
$$   \dot{\bm{p}}_I = \dot{\bm{p}} - \nabla \cdot \bm{\sigma} - \bm{f} , \quad
     \dot{\bm{\pi}}_I = \dot{\bm{\pi}} - \bm{\sigma}
$$

We denote $I_K$ the indicator function of the subset $K$ of a vector space $V$ such that $I_K (\bm{v})$ is equal to zero if $\bm{v} \in K$ and equal to $+\infty$ otherwise. Its polar function is the support function of $K$:
$$ (I_K)^* (\bm{w}) 
     = \sup \left\{ \left\langle \bm{v}, \bm{w} \right\rangle  \mbox{ : } \bm{v} \in K \right\}
$$
In particular $(I_{\left\lbrace \bm{0} \right\rbrace})^* = 0$, then $b_0 (\bm{v}, \bm{w}) = I_{\left\lbrace \bm{0}\right\rbrace} (\bm{v})$ is a  separate bipotential.

Let $E = E^*= S^2(\mathbb{R}^3)$ be the vector space of symmetric tensors of rank 2. 
We consider now a structure of which the material is elastoplastic with a non associated flow rule generated by a bipotential 
$$ b_p : E \times E^* \rightarrow \mathbb{R}\cup \left\{+\infty\right\} : (\bm{e}^p , \bm{s}) \mapsto b_p  (\bm{e}^p , \bm{s})
$$
Let $X = Y$ be the vector space $\mathbb{R}^3 \times S^2(\mathbb{R}^3)$ and $Z$  be the vector space $X \times Y$  of which the elements are denoted  $\dot{\bm{z}} = (\dot{\bm{x}}, \dot{\bm{y}}) = ((\dot{\bm{u}}, \dot{\bm{e}}^p), (\dot{\bm{p}}, \dot{\bm{\pi}})) $. It is straightforward to verify that 
\begin{eqnarray}
    b : Z^2  \rightarrow \mathbb{R}\cup \left\{+\infty\right\} : 
      (\dot{\bm{z}}', \dot{\bm{z}}) \mapsto b (\dot{\bm{z}}', \dot{\bm{z}}) 
\qquad \qquad \quad \nonumber \\
 \quad      =   b_p  (\dot{\bm{\pi}}' , \dot{\bm{\pi}}) + I_{\left\lbrace \bm{0} \right\rbrace} (\dot{\bm{u}}') 
          + I_{\left\lbrace \bm{0} \right\rbrace} ((\dot{\bm{\varepsilon}}^p)')  + I_{\left\lbrace \bm{0} \right\rbrace} (\dot{\bm{p}}') \nonumber
\end{eqnarray}
is a bipotential for the scalar product (\ref{scalar product << . , . >>}). Hence the function $\hat{b}$ given by (\ref{hat(b) (dot(z)_I, dot(z)) = b (J^(-1) dot(z)_I, dot(z)) }) 
\begin{eqnarray}
     \hat{b} : Z^2  \rightarrow \mathbb{R}\cup \left\{+\infty\right\} : 
      (\dot{\bm{z}}_I, \dot{\bm{z}}) \mapsto \hat{b}  (\dot{\bm{z}}_I, \dot{\bm{z}}) \qquad \qquad \quad \nonumber \\ 
      \quad =   b_p  ((\dot{\bm{\varepsilon}}^p)_I , \dot{\bm{\pi}}) + I_{\left\lbrace \bm{0} \right\rbrace} (\dot{\bm{u}}_I) 
          + I_{\left\lbrace \bm{0} \right\rbrace} ( -  \dot{\bm{p}}_I)  + I_{\left\lbrace \bm{0} \right\rbrace} ( - \dot{\bm{\pi}}_I)  \nonumber
\end{eqnarray}
is a symplectic bipotential for the canonical symplectic form (\ref{omega(z, z') = << z, J z' >>}). We are able now to construct the corresponding SBEN principle (\ref{SBEN 2 Pi (z) =}) with
$$ \hat{B}  (\dot{\bm{z}}_I, \dot{\bm{z}}) = \int_\Omega \hat{b}   (\dot{\bm{z}}_I, \dot{\bm{z}})  \,  d\Omega
$$
As the minimum of the functional is the finite value zero, it will be reached when the arguments of the indicator functions vanish. Then we may take the zero value of the indicator functions in the functional  while we introduce extra corresponding constraints
\begin{equation}
    \bm{p} = \rho \, \dot{\bm{u}}, \qquad
 \nabla \cdot \bm{\sigma} +  \bm{f} = \dot{\bm{p}}, \qquad
 \dot{\bm{\pi}} = \bm{\sigma}
 \label{extra constraints}
\end{equation}
Next we can transform the last term in the functional, taking into account the velocity decomposition, the linearity and the antisymmetry of the symplectic form
\begin{eqnarray}
    - \omega (\dot{\bm{z}}_I, \dot{\bm{z}})
 = - \omega (\dot{\bm{z}} - \bm{X}_H, \dot{\bm{z}})
 = - \omega (\dot{\bm{z}}, \dot{\bm{z}})
   + \omega (\bm{X}_H, \dot{\bm{z}}) \nonumber \\
    - \omega (\dot{\bm{z}}_I, \dot{\bm{z}})
 = \omega (\bm{X}_H, \dot{\bm{z}}) \qquad \qquad \quad
 \label{- omega (dot(z)_I, dot(z)) =}
\end{eqnarray}
or in detail
$$ - \omega (\dot{\bm{z}}_I, \dot{\bm{z}})
= \left\langle \bm{p} / \rho, \dot{\bm{p}} \right\rangle
 - \left\langle \dot{\bm{u}}, \nabla \cdot \bm{\sigma} +  \bm{f}  \right\rangle
 - \left\langle \dot{\bm{\varepsilon}}^p, \bm{\sigma} \right\rangle
$$
Owing to the constraints (\ref{extra constraints}), it holds
$$ - \omega (\dot{\bm{z}}_I, \dot{\bm{z}})
= \left\langle \dot{\bm{u}} , \dot{\bm{p}} - \nabla \cdot \bm{\sigma} -  \bm{f}  \right\rangle
 - \left\langle \dot{\bm{\varepsilon}}^p, \bm{\sigma} \right\rangle
= - \left\langle \dot{\bm{\varepsilon}}^p, \bm{\sigma} \right\rangle
$$
Moreover, thanks to the former and latter constraints in (\ref{extra constraints}), the momenta can be eliminated from the functional and the intermediate constraint becomes
$$ \nabla \cdot \bm{\sigma} +  \bm{f} = \rho \, \ddot{\bm{u}} 
$$
while the initial condition on the linear momentum can be transformed into an initial condition on the velocity because $\bm{p} (0) = \rho \, \dot{\bm{u}} (0)$. Taking into account these transformations, the minimum can be searched only on the space of the degrees of freedom $\bm{u}$ and $\bm{\varepsilon}^p $ or equivalently, owing to the elasticity law (\ref{elasticity law}), within the space of couples $(\bm{u}, \bm{\sigma})$. Then we obtain the variational principle:

\textbf{SBEN principle for non associated plasticity.} \textit{The natural evolution of the system minimizes}
\begin{equation}
\Pi(\bm{u}, \bm{\sigma}) = \int_{0}^{T} \left\{ 
\hat{B}  (\nabla \dot{\bm{u}} - \bm{S} : \dot{\bm{\sigma}}, \bm{\sigma}) 
- \left\langle \nabla \dot{\bm{u}} - \bm{S} : \dot{\bm{\sigma}}, \bm{\sigma} \right\rangle
\right\} \mbox{ d}  t
\label{SBEN 3 Pi (x) =}
\end{equation}
\textit{among the admissible paths} $t \mapsto (\bm{u} (t) , \bm{\sigma} (t) ) $ \textit{such that} 
$$ \nabla \cdot \bm{\sigma} +  \bm{f} = \rho \, \ddot{\bm{u}}
$$
\textit{and the minimum is zero.} 

\vspace{0.3cm}

When the bipotential is separated, we recover the variational principle proposed in  \cite{SBEN} and used for numerical simulations 
 in statics for thick tubes \cite{Cao 2020} and plates \cite{Cao 2021a}, and in dynamics \cite{Cao 2021b}.

\vspace{0.5cm}

\subsection{Application to the crack extension with friction contact between the crack sides}

Another important example of non associated flow rule is the model of the unilateral contact with Coulomb's dry friction law. To illustrate it, we consider the problem of extension of a crack in a brittle elastic material in small strain with possible friction contact between the crack sides. In \cite{de Saxce 2022}, we proposed a SBEN principle for the problem of the crack extension without friction introducing as dual variables a crack flow field and a driving force. We recall the main features of this formulation while incorporating the friction contact phenomenon. 

We exclude the event of crack bifurcation into various branches. We denote $S$ the set of admissible surfaces, \textit{i.e.} closed countably $2$ rectifiable subsets of $\Omega$ without change of topology. The crack evolution is described by a time-parameterized family of surfaces
$$ \Gamma : \left[ 0, T \right] \rightarrow S : t \mapsto \Gamma_t = \Gamma (t)
$$
such that $\Gamma(0) = \Gamma_0$ and the map $\Gamma$ is monotone increasing, $S$ being equipped with the inclusion order. The crack front $c_t$ at time $t$ is parameterized by the arc length $s$. The crack extension is modelized by a flow on the cracked surface $\Gamma_T \backslash \Gamma_0$ during the time interval 
\begin{equation}
\left[ 0, T \right] \times c_0 \rightarrow \Gamma_T \backslash \Gamma_0 : (t, \bm{x}_0) \mapsto \bm{x} = \bm{\psi} (t, \bm{x}_0)
\label{crack flow}
\end{equation}
such that 
$$\bm{\psi} (\lbrace t \rbrace \times c_0)= c_t, \qquad
  \bm{\psi} (\left[ 0, t \right] \times c_0)= \Gamma_t \backslash \Gamma_0
$$
The cracked solid at time $t$ is denoted $\Omega_t = \Omega \backslash \Gamma_t$. Within the body, the displacement field  at time $t$ is 
$$ \Omega_t  \rightarrow \mathbb{R}^3 : \bm{x} \mapsto \bm{u} (t, \bm{x})
$$
A crack being a material discontinuity, we must distinguish  two material surfaces $\Gamma^+_t$ and $\Gamma^-_t$ that occupy the same position as $\Gamma_t$ but are the two sides of the crack and have opposite unit normal vectors, exterior to $\Omega_t$: $\bm{n} = \bm{n}^- = - \bm{n}^+$. Let $\bm{u}^\pm$ (resp. $\bm{\sigma}^\pm $) be  the value of the displacement (resp. stress tensor) on $\Gamma^\pm_t$, $\left[\bm{u}\right] = \bm{u}^+ - \bm{u}^-$ be the relative displacement and $\bm{t} = \bm{\sigma}^+ \cdot \bm{n}$ be the reaction subjected to $\Gamma^+_t$ from $\Gamma^-_t$.  The vector $\left[\bm{u}\right] $ (resp. $ \bm{t}$) can be decomposed into its tangential part $\left[\bm{u}_t \right] $ (resp. $ \bm{t}_t$) and its normal part $\left[u_n \right] $ (resp. $t_n$). If the two sides are in contact, Coulomb's law states that
\begin{eqnarray}
 \mbox{if} \;\bm{t} = 0 & , & \quad   \left[\dot{\bm{u}}\right]  > 0,\nonumber\\
      \mbox{if} \;  \parallel \bm{t}_t \parallel < \mu \, t_n & , & \quad  \left[\dot{\bm{u}}\right]  = \bm{0}   \nonumber\\
      \mbox{if} \;  \parallel \bm{t}_t \parallel = \mu \, t_n & , & \quad  \exists \lambda \geq 0 \; \mbox{such that} \;
                 - \left[\dot{\bm{u}}\right]  = \lambda \, \frac{\bm{t}_t}{\parallel  \bm{t}_t \parallel}\nonumber
\end{eqnarray}
In (\cite{saxfeng}, \cite{sax CRAS 92}, \cite{sax boussh IJMS 98}),  we proved that this law can be generated by the bipotential
$$  b_c : \mathbb{R} ^3  \times  \mathbb{R} ^3 \rightarrow \mathbb{R}\cup \left\{+\infty\right\}
$$
which has a finite value
$$ b_c ( - \left[\dot{\bm{u}}\right]  , \bm{t}) = \mu \, t_n \parallel  \left[\dot{\bm{u}}_t \right]  \parallel
$$
provided $ \parallel \bm{t}_t \parallel \leq \mu \, t_n$ and $  \left[ u_n \right]   \geq 0$.

Moreover, by comparison to the experimental data and the crack stability criteria proposed in the literature, we showed in \cite{de Saxce 2022} the relevance to represent the crack extension in brittle materials by a normality law of the form 
$$ \dot{\bm{\psi}} \in \partial \varphi (\bm{G})
$$
where $\varphi$ is a convex lower semicontinuous function, in fact the indicator function of the closed convex domain of crack stability.

The elements of the space $X = Y = \mathbb{R}^3 \times \mathbb{R}^3 \times \mathbb{R}^3$ are the triplets $\bm{\xi} =  (\bm{u}, - \left[\bm{u}\right] , \bm{\psi}) $ of fields of displacement, relative displacement and crack flow. The corresponding dynamical momenta are denoted $\bm{\eta} = (\bm{p} , \bm{\pi}_c, \bm{\pi}_f) $. The dual pairing is of the form
$$ \left\langle \dot{\bm{\xi}} , \dot{\bm{\eta}} \right\rangle 
 = \left\langle \dot{\bm{u}} , \dot{\bm{p}} \right\rangle + \left\langle - \left[\dot{\bm{u}}\right]  , \dot{\bm{\pi}_c} \right\rangle 
     + \left\langle \dot{\bm{\psi}} , \dot{\bm{\pi}}_f \right\rangle 
$$
$$ \left\langle \dot{\bm{\xi}} , \dot{\bm{\eta}} \right\rangle 
= \int_{\Omega_t} \dot{\bm{u}} \cdot \dot{\bm{p}}  \, d\Omega
  - \int_{\Gamma_t}  \left[\dot{\bm{u}}\right]  \cdot \dot{\bm{\pi}}_c  \, \mbox{d} \Gamma
  + \int_{c_t} \dot{\bm{\psi}} \cdot \dot{\bm{\pi}}_f \, \mbox{d} s
$$
The Hamiltonian of the structure is 
\begin{eqnarray}
    H (t, \bm{z}) 
    = \int_{\Omega} \left\lbrace \dfrac{1}{2 \rho} \parallel \bm{p} \parallel^ 2 
    + w\, (\nabla \bm{u} - \bm{\varepsilon}^p) - \bm{f} (t)\cdot\bm{u} \right\rbrace  \, d\Omega \nonumber \\
              - \int_{\Gamma_1} \bar{\bm{f}} (t)\cdot \bm{u}  \, d\Gamma \nonumber
\end{eqnarray}
We consider the canonical symplectic form
\begin{equation}
   \omega \, (\dot{\bm{z}}, \dot{\bm{z}}') =  \omega \, ( (\dot{\bm{\xi}}, \dot{\bm{\eta}}) , (\dot{\bm{\xi}}', \dot{\bm{\eta}}'))
    = \left\langle  \dot{\bm{\xi}}, \dot{\bm{\eta}}'  \right\rangle - \left\langle  \dot{\bm{\xi}}', \dot{\bm{\eta}}  \right\rangle
\label{omega(dot(z), dot(z)') = < xi , eta' > - < xi' , eta >}
\end{equation}
The symplectic gradient of the Hamiltonian is
$$ \bm{X}_H = ((D_{\bm{p}} H, D_{\bm{\pi}_c} H, D_{\bm{\pi}_f} H), 
                         (- D_{\bm{u}} H, - D_{- \left[\bm{u}\right] } H, D_{\bm{\psi}} H ))
$$
$$ \bm{X}_H = ( ( \bm{p} / \rho, \bm{0} , \bm{0}) , ( \nabla \cdot \bm{u} + \bm{f}, \bm{t} , \bm{G} ) )
$$
where $\bm{G}$ is the driving force of which the explicit expression was calculated in \cite{de Saxce 2022}.
Thus one has
$$ \dot{\bm{z}}_I = \dot{\bm{z}} - \bm{X}_H 
                            = ( ( \dot{\bm{u}}_I , - \left[\dot{\bm{u}}\right] _I,  \dot{\bm{\psi}}_I  ),     
                                 ( \dot{\bm{p}}_I, (\dot{\bm{\pi}}_c)_I, (\dot{\bm{\pi}}_f)_I ) )
$$
with
$$ \dot{\bm{u}}_I =\dot{\bm{u}} -  \bm{p} / \rho , \quad
     - \left[\dot{\bm{u}}\right] _I =  - \left[\dot{\bm{u}}\right] , \quad 
     \dot{\bm{\psi}}_I =\dot{\bm{\psi}} ,
$$
$$     \dot{\bm{p}}_I = \dot{\bm{p}} - \nabla \cdot \bm{u} - \bm{f} , \quad
     (\dot{\bm{\pi}}_c)_I=\dot{\bm{\pi}}_c - \bm{t}  , \quad
     (\dot{\bm{\pi}}_f)_I = \dot{\bm{\pi}}_f - \bm{G}
$$
Let $Z$  be the vector space $X \times Y$  of which the elements are denoted  $\dot{\bm{z}} = (\dot{\bm{\xi}}, \dot{\bm{\eta}}) = (  (\bm{u}, - \left[\bm{u}\right] , \bm{\psi}) , (\bm{p} , \bm{\pi}_c, \bm{\pi}_f)  ) $. We can readily verify that 
\begin{eqnarray}
   b : & & Z^2  \rightarrow \mathbb{R}\cup \left\{+\infty\right\} : 
      (\dot{\bm{z}}', \dot{\bm{z}}) \mapsto b (\dot{\bm{z}}', \dot{\bm{z}}) \nonumber\\ 
    & =  & b_c  (\dot{\bm{\pi}}'_c , \dot{\bm{\pi}}_c ) 
          + \varphi (\dot{\bm{\pi}}_f) + \varphi^* (\dot{\bm{\pi}}'_f) \nonumber\\
    &   &  + I_{\left\lbrace \bm{0} \right\rbrace} (\dot{\bm{u}}') 
          + I_{\left\lbrace \bm{0} \right\rbrace} (- \left[\dot{\bm{u}}' \right])
          + I_{\left\lbrace \bm{0} \right\rbrace} ((\dot{\bm{\psi}})')  
          + I_{\left\lbrace \bm{0} \right\rbrace} (\dot{\bm{p}}') \nonumber
\end{eqnarray}
is a bipotential for the scalar product (\ref{scalar product << . , . >>}). Using (\ref{hat(b) (dot(z)_I, dot(z)) = b (J^(-1) dot(z)_I, dot(z)) }),  the function
\begin{eqnarray}
   \hat{b} : & & Z^2  \rightarrow \mathbb{R}\cup \left\{+\infty\right\} : 
      (\dot{\bm{z}}', \dot{\bm{z}}) \mapsto b (\dot{\bm{z}}', \dot{\bm{z}}) \nonumber\\ 
    & =  & b_c  (  - \left[\dot{\bm{u}} \right]_I , \dot{\bm{\pi}}_c ) 
          + \varphi (\dot{\bm{\pi}}_f) + \varphi^* (\dot{\bm{\psi}}_I) \nonumber\\
   &  &   + I_{\left\lbrace \bm{0} \right\rbrace} (\dot{\bm{u}}_I) 
          + I_{\left\lbrace \bm{0} \right\rbrace} (- \dot{\bm{p}}_I) \nonumber\\
   &  &   + I_{\left\lbrace \bm{0} \right\rbrace} ( -(\dot{\bm{\pi}_c})_I)  
          + I_{\left\lbrace \bm{0} \right\rbrace} ( -(\dot{\bm{\pi}_f})_I)    \nonumber
\end{eqnarray}
is a symplectic bipotential for the canonical symplectic form (\ref{omega(dot(z), dot(z)') = < xi , eta' > - < xi' , eta >}). 

Next, we take the zero value of the indicator functions in the functional  while we introduce extra corresponding constraints
\begin{equation}
    \bm{p} = \rho \, \dot{\bm{u}}, \qquad
 \nabla \cdot \bm{\sigma} +  \bm{f} = \dot{\bm{p}}, \qquad
 \dot{\bm{\pi}}_c = \bm{t},  \qquad
  \dot{\bm{\pi}}_f = \bm{G}
 \label{extra constraints 2}
\end{equation}
Similarly to the calculations of Section 'Application to the non associated plasticity', we can transform the last term in the functional, taking into account the previous constraints
$$    - \omega (\dot{\bm{z}}_I, \dot{\bm{z}})
 = \omega (\bm{X}_H, \dot{\bm{z}})
= \left\langle \left[ \dot{\bm{u}} \right] , \bm{t}  \right\rangle
    - \left\langle \dot{\bm{\psi}}, \bm{G} \right\rangle
$$
Finally, we derive the corresponding SBEN principle from (\ref{SBEN 2 Pi (z) =})

\textbf{SBEN principle for crack extension and friction contact.} \textit{The natural evolution of the system minimizes}
\begin{eqnarray}
\Pi(\bm{u}, \bm{\psi}) = \int_{0}^{T} \{ 
\int_{\Gamma_t} (b_c (- \left[ \dot{\bm{u}} \right] , \bm{t} ) + \left[ \dot{\bm{u}} \right] \cdot \bm{t} ) \, \mbox{d} \Gamma \nonumber \\
  + \int_{c_t} ( \varphi (\bm{G}) + \varphi^* (\dot{\bm{\psi}}) - \dot{\bm{\psi}} \cdot \bm{G}) \, \mbox{d} s 
\} \mbox{ d}  t 
\label{SBEN 4 Pi (x) =}
\end{eqnarray}
\textit{among the admissible paths} $t \mapsto (\bm{u} (t) , \bm{\sigma} (t) ) $
\textit{such that} 
$$  \nabla \cdot \bm{\sigma} +  \bm{f} = \rho \, \ddot{\bm{u}}
$$
\textit{and the minimum is zero.} 

\section{Conclusions}

\label{Section - Conclusions}

In previous papers, we used the Brezis-Ekeland-Nayroles principle (BEN) for numerical applications to the standard plasticity and viscoplasticity with the normality law in dynamics thanks to a symplectic version SBEN proposed in \cite{Buliga 2008}.

In this work, we generalized this variational principle to the non associated plasticity using the concept of bipotential. Its spectrum of applications to constitutive laws is much broader than the one of maximal monotone laws for which Visintin extended the BEN principle  (\cite{Visintin 2008}, \cite{Visintin 2013}).

For this aim, we introduced the general concept of symplectic bipotential and specialized it to two kinds of applications, the non associated plasticity and the crack extension with friction contact between the crack sides. In the future, we  intend to develop for these two applications numerical simulations based on this new variational principle. 

\vspace{0.5cm}

\textbf{Aknowledgments}

\vspace{0.3cm}

This work was performed thanks to the research  project BIpotentiels G\'{e}n\'{e}ralis\'{e}s pour le principe variationnel de Brezis-Ekeland-Nayroles en m\'{e}canique (BigBen) supported by the Agence Nationale de la Recherche (Project ANR-22-CE51-0034-04).


\end{document}